\documentclass[11pt]{article}
\usepackage{amsfonts}
\usepackage{amsmath}

\newtheorem{theorem}{Theorem}

\newtheorem{corollary}[theorem]{Corollary}

\input{tcilatex}

\begin{document}

\title{Singular elliptic equation involving the GJMS operator on the
standard unit sphere.}
\author{Mohammed Benalili and Ali Zouaoui}
\maketitle

\begin{abstract}
Given a Riemannian compact manifold $\left( M,g\right) $ of dimension $n\geq
5$, we have proven in \cite{1} under some conditions that the equation : 
\begin{equation}
P_{g}(u)=Bu^{2^{\sharp }-1}+\frac{A}{u^{2^{\sharp }+1}}+\frac{C}{u^{p}}
\label{E1}
\end{equation}%
where $P_{g}$ is GJMS-operator, $n=\dim (M)>2k$ $(k\in \mathbb{N}^{\star })$%
, $A,B$ and $C$ are smooth positive functions on $M$, $p>1$ and $2^{\sharp }=%
\frac{2n}{n-2k}$ denotes the critical Sobolev of the embedding $H_{k}^{2}(M)$
$\subset $ $L^{2^{\sharp }}(M)$, admits two distinct positive solutions. The
proof of this result is essentially based on the given smooth function $%
\varphi >0$ with norm $\Vert \varphi \Vert _{P_{g}}=1$ fulfilling some
conditions ( see Theorem 3 in \cite{1}). In this note we construct an
example of such function on the unit standard sphere $\left( \mathbb{S}%
^{n},h\right) $. Consequently the conditions of the Theorem are improved in
the case of $\left( \mathbb{S}^{n},h\right) $.
\end{abstract}

\section{Construction of the function $\protect\varphi $}

Inspired by the work of F. Robert. (see \cite{4} ) , we construct an example
of a smooth function $\varphi >0$ on the Euclidean sphere $\left( \mathbb{S}%
^{n},h\right) $ with norm$\Vert \varphi \Vert _{P_{h}}=1$.\newline
Indeed let $\lambda >0$ and $x_{0}\in \mathbb{S}^{n}$. To a rotation, we may
assume that $x_{0}$ is the north pole i.e. $x_{0}=\left( 0,...0,1\right) $.
We consider the transformation 
\begin{equation*}
\phi _{\lambda }:\mathbb{S}^{n}\rightarrow \mathbb{S}^{n}
\end{equation*}%
defined by $\phi _{\lambda }(x)=\psi _{x_{0}}^{-1}\left( \lambda ^{-1}.\psi
_{x_{0}}(x)\right) $ if $x\neq x_{0}$ and $\phi _{\lambda }(x_{0})=x_{0}$
where $\psi _{x_{0}}$ is the stereographic projection of $x_{0}$ given by 
\begin{equation*}
\psi _{x_{0}}:\left( \mathbb{S}^{n}\setminus \left\{ x_{0}\right\} ,h\right)
\rightarrow \left( \mathbb{R}^{n},\xi \right) ,
\end{equation*}%
for any $a=\left( \eta _{1},...,\eta _{n},\zeta \right) $ associates $\psi
_{x_{0}}(a)=\left( \frac{\eta _{1}}{1-\zeta },...,\frac{\eta _{n}}{1-\zeta }%
\right) $ and 
\begin{equation*}
\begin{array}{c}
\delta _{\lambda }:\left( \mathbb{R}^{n},\xi \right) \rightarrow \left( 
\mathbb{R}^{n},\xi \right)  \\ 
x\mapsto \delta _{\lambda }(x)=\frac{1}{\lambda }x%
\end{array}%
\end{equation*}%
is the homothetic mapping. $h$ is the canonical metric on $\mathbb{S}^{n}$
and $\xi $ is the Euclidean one on $\mathbb{R}^{n}$.\newline
\newline
Note that $\psi _{x_{0}}$ is a conformal, mapping more precisely we have 
\begin{equation*}
\left( \psi _{x_{0}}^{-1}\right) ^{\star }h=U^{\frac{4}{n-2k}}.\xi 
\end{equation*}%
where $U(x)=\left( \frac{1+\Vert x\Vert ^{2}}{2}\right) ^{k-\frac{n}{2}}$ .
Hence $\phi _{\lambda }$ is conformal i.e. 
\begin{equation*}
\phi _{\lambda }^{\star }h=u_{x_{0},\beta }^{\frac{4}{n-2k}}.h\quad \text{%
where}\;\beta =\frac{1+\lambda ^{2}}{\lambda ^{2}-1}
\end{equation*}%
and 
\begin{equation*}
u_{x_{0},\beta }(x)=\left( \dfrac{\sqrt{\beta ^{2}-1}}{\beta -\cos
d_{h}(x,x_{0})}\right) ^{\frac{n-2k}{2}}\quad \forall x\in \mathbb{S}^{n}\;%
\text{with}\;\beta >1.
\end{equation*}%
In particular we have 
\begin{equation}
\int\limits_{\mathbb{S}^{n}}u_{x_{0},\beta }^{2^{\sharp }}dv_{h}=\omega _{n}
\label{01}
\end{equation}%
where $\omega _{n}>0$ is the volume of the unit standard sphere $\left( 
\mathbb{S}^{n},h\right) .$\newline
By the conformal invariance of the operator $P_{h}$ on $\left( \mathbb{S}%
^{n},h\right) $, we obtain that 
\begin{equation}
P_{h}(u_{x_{0},\beta })=\frac{n-2k}{2}Q_{h}u_{x_{0},\beta }^{2^{\sharp }-1}
\label{02}
\end{equation}%
where $Q_{h}$ denotes the $Q$-curvature of $\left( \mathbb{S}^{n},h\right) $
which expresses by the Gover's formula as: 
\begin{equation*}
Q_{h}=\dfrac{2}{n-2k}P_{h}(1)=\dfrac{2}{n-2k}(-1)^{k}\prod_{l=1}^{k}(c_{l}%
\;Sc)
\end{equation*}%
where $c_{l}=\frac{(n+2l-2)(n-2l)}{4n(n-1)}$, $Sc=n(n-1)$ (the scalar
curvature of $\left( \mathbb{S}^{n},h\right) $). So the $Q_{h\text{ \ }}$is
a positive constant.\newline
Multiplying the two sides of \eqref{02} by $u_{x_{0},\beta }$ and
integrating on $\mathbb{S}^{n}$ we get: 
\begin{equation*}
\int\limits_{\mathbb{S}^{n}}u_{x_{0},\beta }P_{h}(u_{x_{0},\beta })dv_{h}=%
\frac{n-2k}{2}Q_{h}\int\limits_{\mathbb{S}^{n}}u_{x_{0},\beta }^{2^{\sharp
}}dv_{h}.
\end{equation*}%
And since 
\begin{equation*}
\int\limits_{\mathbb{S}^{n}}u_{x_{0},\beta }P_{h}(u_{x_{0},\beta
})dv_{h}=\Vert u_{x_{0},\beta }\Vert _{P_{h}}^{2}
\end{equation*}%
\eqref{02} writes 
\begin{equation*}
\Vert u_{x_{0},\beta }\Vert _{P_{h}}^{2}=\frac{n-2k}{2}Q_{h}\omega _{n}.
\end{equation*}%
Hence by putting 
\begin{equation*}
\varphi =\left( \frac{n-2k}{2}\omega _{n}Q_{h}\right) ^{\frac{-1}{2}%
}u_{x_{0},\beta }
\end{equation*}%
we obtain a function satisfying the conditions of Theorem 3 in \cite{1} i.e. 
$\varphi >0$ smooth on $(\mathbb{S}^{n},h)$ such that $\Vert \varphi \Vert
_{P_{g}}=1$

\section{ Existence results on the sphere}

On the standard unit sphere $(\mathbb{S}^{n},h),$ if we take the function $%
\varphi $ of Theorem 3 in \cite{1} equals $\left( \frac{n-2k}{2}\omega
_{n}Q_{h}\right) ^{\frac{-1}{2}}u_{x_{0},\beta }$ , we obtain

\begin{theorem}
\label{th1}{\ }Let $\left( \mathbb{S}^{n},h\right) $ be the unit standard
unit sphere of dimension $n>2k,\;k\in \mathbb{N}^{\star }.$ There is a
constant $C(n,p,k)>0$ depending only on $n,p,k$ such that 
\begin{equation}
\frac{1}{2^{\sharp }}\left( \frac{n-2k}{2}\omega _{n}Q_{h}\right) ^{\frac{%
2^{\sharp }}{2}}\int_{\mathbb{S}^{n}}\frac{A(x)}{u_{x_{0},\beta
}^{2^{\natural }}}dv_{h}\leq C\left( n,p,k\right) \left( S\underset{x\in 
\mathbb{S}^{n}}{\max }B(x)\right) ^{\frac{2+2^{\sharp }}{2-2^{\sharp }}}
\label{2.3}
\end{equation}%
and 
\begin{equation}
\frac{1}{p-1}\left( \frac{n-2k}{2}.\omega _{n}.Q_{h}\right) ^{\frac{p-1}{2}%
}\int_{\mathbb{S}^{n}}\frac{C(x)}{u_{x_{0},\beta }^{p-1}}dv_{h}\leq C\left(
n,p,k\right) \left( S\underset{x\in \mathbb{S}^{n}}{\max }B(x)\right) ^{%
\frac{p+1}{2-2^{\sharp }}}  \label{2.4}
\end{equation}%
where 
\begin{equation*}
u_{x_{0},\beta }(x)=\left( \dfrac{\sqrt{\beta ^{2}-1}}{\beta -\cos
d_{h}(x,x_{0})}\right) ^{\frac{n-2k}{2}}\quad \forall x\in \mathbb{S}^{n}\;%
\text{and}\;\beta >1.
\end{equation*}%
then the equation \eqref{E1} admits a solution of class $C^{\infty }(\mathbb{%
S}^{n})$. If moreover for any $\varepsilon \in \left] 0,\lambda ^{\star }%
\right[ $ where $\lambda ^{\star }$ is a positive constant the two following
conditions are satisfied 
\begin{equation*}
\frac{2}{a}\left( \int\limits_{\mathbb{S}^{n}}\sqrt{A(x)}dv_{h}\right)
^{2}\left( \frac{1}{t_{0}a_{1}}\right) ^{2^{\sharp }}>2^{\sharp }k\frac{%
t_{0}^{2}}{4n}(2-a)
\end{equation*}%
and 
\begin{equation*}
\left( \frac{2}{a}\right) ^{\frac{p-1}{2^{\sharp }}}\left( \int\limits_{%
\mathbb{S}^{n}}\sqrt{C(x)}dv_{h}\right) ^{2}\left( \frac{1}{t_{0}a_{2}}%
\right) ^{p-1}>(p-1)k\frac{t_{0}^{2}}{4n}(2-a)
\end{equation*}%
where $a_{1},a_{2}$ are positive constants, $2^{\sharp }=\frac{2n}{n-2k}%
\;,3<p<2^{\sharp }+1$. Then the equation \eqref{E1} admits a second solution.
\end{theorem}

Note that since 
\begin{equation}
\left( \dfrac{\beta -1}{\beta +1}\right) ^{\frac{n-2k}{4}}\leq
u_{x_{0},\beta }(x)\leq \left( \dfrac{\beta +1}{\beta -1}\right) ^{\frac{n-2k%
}{4}}  \label{03}
\end{equation}%
we can improve the conditions (\ref{2.3}) and (\ref{2.4}) of Theorem \ref%
{th1}. Indeed, from (\ref{0.3}) we deduce that 
\begin{equation*}
\varphi (x)\geq \left( \dfrac{\beta -1}{\beta +1}\right) ^{\frac{n-2k}{4}%
}\left( \frac{n-2k}{2}\omega _{n}Q_{h}\right) ^{-\frac{1}{2}}.
\end{equation*}%
Consequently 
\begin{equation*}
\frac{\Vert \varphi \Vert ^{2^{\sharp }}}{2^{\sharp }}\int_{\mathbb{S}^{n}}%
\frac{A(x)}{\varphi ^{2^{\natural }}}dv_{h}=\frac{1}{2^{\sharp }}\int_{%
\mathbb{S}^{n}}\frac{A(x)}{\varphi ^{2^{\natural }}}dv_{h}\leq \frac{1}{%
2^{\sharp }}\left( \dfrac{\beta +1}{\beta -1}\right) ^{\frac{n}{2}}\left( 
\frac{n-2k}{2}.\omega _{n}.Q_{h}\right) ^{\frac{2^{\sharp }}{2}}\int_{%
\mathbb{S}^{n}}A(x)dv_{h}
\end{equation*}%
So, if 
\begin{equation*}
\frac{1}{2^{\sharp }}\left( \dfrac{\beta +1}{\beta -1}\right) ^{\frac{n}{2}%
}\left( \frac{n-2k}{2}.\omega _{n}.Q_{h}\right) ^{\frac{2^{\sharp }}{2}%
}\int_{\mathbb{S}^{n}}A(x)dv_{h}\leq C\left( n,p,k\right) \left( S\underset{%
x\in \mathbb{S}^{n}}{\max }B(x)\right) ^{\frac{2+2^{\sharp }}{2-2^{\sharp }}}
\end{equation*}%
then the condition (\ref{2.3}) is fulfilled. Likewise if 
\begin{equation*}
\frac{1}{p-1}\left( \frac{n-2k}{2}\omega _{n}Q_{h}\right) ^{\frac{p-1}{2}%
}\left( \dfrac{\beta +1}{\beta -1}\right) ^{\frac{n(p-1)}{2.2^{\sharp }}%
}\int_{\mathbb{S}^{n}}C(x)dv_{h}\leq C\left( n,p,k\right) \left( S\underset{%
x\in M}{\max }B(x)\right) ^{\frac{p+1}{2-2^{\sharp }}}
\end{equation*}%
the condition (\ref{2.4}) is also true and we deduce the following result:

\begin{corollary}
{\ }Let $\left( \mathbb{S}^{n},h\right) $ be the unit standard unit sphere
of dimension $n>2k,\;k\in \mathbb{N}^{\star }.$ There is a constant $%
C(n,p,k)>0$ depending only on $n,p,k$ such that 
\begin{equation}
\frac{1}{2^{\sharp }}\left( \dfrac{\beta +1}{\beta -1}\right) ^{\frac{n}{2}%
}\left( \frac{n-2k}{2}.\omega _{n}.Q_{h}\right) ^{\frac{2^{\sharp }}{2}%
}\int_{\mathbb{S}^{n}}A(x)dv_{h}\leq C\left( n,p,k\right) \left( S\underset{%
x\in \mathbb{S}^{n}}{\max }B(x)\right) ^{\frac{2+2^{\sharp }}{2-2^{\sharp }}}
\end{equation}%
and 
\begin{equation}
\frac{1}{p-1}\left( \frac{n-2k}{2}.\omega _{n}.Q_{h}\right) ^{\frac{p-1}{2}%
}\left( \dfrac{\beta +1}{\beta -1}\right) ^{\frac{n(p-1)}{2.2^{\sharp }}%
}\int_{\mathbb{S}^{n}}C(x)dv_{h}\leq C\left( n,p,k\right) \left( S\underset{%
x\in M}{\max }B(x)\right) ^{\frac{p+1}{2-2^{\sharp }}}
\end{equation}%
where $\beta >1$. Then the equation\eqref{E1} admits a solution of class $%
C^{\infty }(\mathbb{S}^{n})$. If moreover for any $\varepsilon \in \left]
0,\lambda ^{\star }\right[ $, where $\lambda ^{\star }$ is a positive
constant, the two following assumptions are satisfied 
\begin{equation}
\frac{2}{a}\left( \int\limits_{\mathbb{S}^{n}}\sqrt{A(x)}dv_{h}\right)
^{2}\left( \frac{1}{t_{0}a_{1}}\right) ^{2^{\sharp }}>2^{\sharp }k\frac{%
t_{0}^{2}}{4n}(2-a)
\end{equation}%
and 
\begin{equation}
\left( \frac{2}{a}\right) ^{\frac{p-1}{2^{\sharp }}}\left( \int\limits_{%
\mathbb{S}^{n}}\sqrt{C(x)}dv_{h}\right) ^{2}\left( \frac{1}{t_{0}a_{2}}%
\right) ^{p-1}>(p-1)k\frac{t_{0}^{2}}{4n}(2-a)
\end{equation}%
where $a_{1},a_{2}$ are positive constants, $2^{\sharp }=\frac{2n}{n-2k}%
\;,3<p<2^{\sharp }+1$. Then the equation \eqref{E1} admits a second solution.
\end{corollary}

\end{document}